\documentclass[12pt]{article}

\usepackage{multirow}
\usepackage{graphicx}
\usepackage{amstext}
\usepackage{amsmath}
\usepackage{amssymb}

\usepackage{makecell}
\usepackage{natbib}
\usepackage{url}

\def\proof{\par\smallskip\noindent{Proof.} }
\def\endproof{\hfill $\Box$ \vskip .5cm}

\bibpunct[, ]{[}{]}{,}{a}{}{,}%
\def\bibfont{\small}%
\def\bibsep{\smallskipamount}%
\def\bibhang{24pt}%

\newtheorem{thm}{Theorem}%[section]
\newtheorem{cor}{Corollary}%[section]
\newtheorem{exam}{Example}%[section]

\def\dg{{\rm diag}}

\def\ovl{\overline}

\def\non{\nonumber}

\begin{document}

\begin{center}

{\LARGE  Tight MIQP Reformulations for Semi-Continuous Quadratic Programming:
\\Lift-and-Convexification Approach}\\[12pt]

\footnotesize
\mbox{\large Baiyi Wu}\\
Department of Systems Engineering and Engineering Management, The
Chinese University of Hong Kong, Shatin, N. T., Hong Kong,
\mbox{bywu@se.cuhk.edu.hk} \\
[6pt]

\mbox{\large Xiaoling Sun}\\
Department of Management Science,
School of Management, Fudan University, Shanghai 200433, P. R. China\\
[6pt]

\mbox{\large Duan Li}\\
Department of Systems Engineering and Engineering Management, The
Chinese University of Hong Kong,  Shatin, N. T., Hong Kong,
\mbox{dli@se.cuhk.edu.hk} \\
[6pt]

\mbox{\large Xiaojin Zheng}\\
School of Economics and Management,
Tongji University, Shanghai 200092, P. R. China.\\
[6pt]

\normalsize

\end{center}

%\maketitle
%\slugger{mms}{xxxx}{xx}{x}{x--x}
%slugger should be set to mms, siap, sicomp, sicon, sidma, sima, simax, sinum, siopt, sisc, or sirev

\begin{abstract}
We consider in this paper a class of semi-continuous quadratic programming problems which arises in many real-world applications such as
production planning, portfolio selection and subset selection in regression.
We propose a lift-and-convexification approach
to derive an equivalent reformulation of the original problem.
This lift-and-convexification approach lifts the quadratic term involving $x$
only in the original objective function $f(x,y)$ to a quadratic function
of both $x$ and $y$ and convexifies this equivalent objective function.
While the continuous relaxation of our new reformulation
attains the same tight bound as achieved by the continuous relaxation of the
well known perspective reformulation, the new reformulation also
retains the linearly constrained quadratic programming structure of the original mix-integer problem.
This prominent feature improves the performance of branch-and-bound algorithms by providing the same
tightness at the root node as the state-of-the-art perspective reformulation and offering
much faster processing time at children nodes.
We further combine the lift-and-convexification approach and the quadratic convex reformulation
approach in the literature to form an even tighter reformulation.
Promising results from our computational tests in both portfolio selection and subset
selection problems numerically verify the benefits from these theoretical features
of our new reformulations.

\end{abstract}

%\begin{keywords}
%mixed integer quadratic programming;
%semi-continuous quadratic program;
%perspective cut reformulation;
%quadratic convex reformulation;
%semidefinite program;
%portfolio selection;
%subset selection.
%\end{keywords}

%\begin{AMS}
%90C11, 90C20, 90C34, 90B50.
%\end{AMS}

%\pagestyle{myheadings}
%\thispagestyle{plain}
%\markboth{B. Y. WU, X. L. SUN, D. LI AND X. J. ZHENG}
%{LIFT-AND-CONVEXIFICATION FOR SEMI-CONTINUOUS QP}

\section{Introduction}
We consider in this paper the following mixed-integer quadratic
programming (MIQP) problem:
\begin{align}
({\rm P})~~\non
\min
&~~f(x,y)=x^TQx+c^Tx+h^Ty\non\\
{\rm s.t.}
&~~Ax+By\le d,\non\\
&~~a_iy_i\le x_i\le b_i y_i,~y_i\in\{0,1\},~i=1,\ldots,n, \label{semi}
\end{align}
where $Q$ is an $n\times n$ positive semidefinite symmetric matrix,
$c,h\in\Re^n$, $d\in\Re^m$, and $A,B\in\Re^{m\times n}$.

Problem (P) is in general NP-hard (see \cite{bienstock96}).
Its difficulty arises from the discrete structure induced by
the constraint in (\ref{semi}).
This constraint is used to model the situation where $x_i$ must rest inside
an interval if it is not zero, that is, $x_i\in\{0\}\cup [a_i,b_i]$.
These variables $x_i$ are termed {\em semi-continuous} variables.
We assume in our study $a_i<b_i,~i=1,...,n$.
We also assume that the feasible region of problem (P) is nonempty.
A recent review on problem (P) and its solution methods can be
found in \cite{sun2013}.

Semi-continuous variables appear in many real-world optimization problems.
For instance, in production planning, the semi-continuous variables are used to
describe the state of a production process that is either turned off
(inactive), hence nothing is produced, or turned on (active) such
that the production level has to lie in certain interval
(\cite{frangioni2006or,frangioni2008power,frangioni2009ieee}).
Other typical applications of semi-continuous variables include portfolio
selection with {\em minimum buy-in threshold}
(\cite{jobst01,frangioni06mp,cui2013,sun2013}) and lot-sizing with
minimum order quantity (\cite{anderson1993,park2013lot}).

An important instance of (P) involves optimization models with
a cardinality constraint:
\begin{align}
   |{\rm supp}(x)|\le K,
   \label{card}
\end{align}
where ${\rm supp}(x)=\{i~|~x_i\neq 0\}$ and $K$ is an integer with
$1\le K\le n$. The cardinality constraint is often encountered when
the number of nonzero variables has to be  limited. The cardinality
constraint in (\ref{card}) can be easily incorporated into problem
(P) by introducing an additional linear constraint
$\sum_{i=1}^my_i\le K$.

A well-known application of semi-continuous variables and cardinality
constraint is the {\it cardinality constrained mean-variance portfolio selection} in financial optimization.
The classical mean-variance model of Markowitz is a quadratic programming
problem that minimizes the variance subject to linear
constraints on expected return and budget availabilities.
In real-world applications of portfolio selection models,
however,  most investors would invest in only a limited number of assets
due to market frictions such as management and transaction fees.
Moreover, the minimum buy-in threshold is often a mandate trading constraint.
Suppose that there are $n$ risky assets in a financial market with a
random return vector $R=(R_1,\ldots,R_n)^T$. Furthermore, the expected return
vector and the covariance matrix of $R$ are assumed to be given as $\mu$ and $Q$,
respectively.
The portfolio selection model with cardinality and
minimum buy-in threshold constraints can be then expressed as:
\begin{align*}
{\rm (MV)}~~
\min
&~~x^TQx \\
{\rm s.t.}
&~~\sum^n_{i=1}x_i=1,\\
&~~\mu^Tx\ge \rho,\\
&~~\sum_{i=1}^ny_i\le K,\\
&~~a_iy_i\le x_i\le b_i y_i,~y_i\in\{0,1\},~i=1,\ldots,n,
\end{align*}
where $x_i$ represents the proportion of the total capital invested
in the $i$th asset, and $\rho$ is a prescribed expected return level set by the
investor. Portfolio selection problems with cardinality and/or
minimum threshold constraints have been studied extensively in recent literature.
For exact solution methods, please see, e.g.,
\cite{bienstock96,li06,shawa08,bonami09,bertsimas09,cui2013,gaoli2013b}.
For inexact solution methods, such as heuristics, local search
methods and randomized techniques, please see, e.g.,
\cite{jacob74,blog83,chang00,jobst01,schaerf02,maringer03,crama03,mitra07,fernandez07,zhang08}.

Another application of (P) with cardinality constraint is the subset selection problem
in multivariate linear regression.
Given $m$ observed data points $(a_i,b_i)$ with $a_i\in\Re^n$ and $b_i\in\Re$, we need
to minimize the least square measure of $\sum_{i=1}^{m}(a_i^Tx-b_i)^2$ with only
a subset of the prediction variables in $x$
(see, e.g., \cite{arthanari93,miller02,bertsimas09}).
This problem can be formally formulated as:
\begin{align*}
{\rm (SSP)}~~
\min
&~~\|Ax-b\|^2 \\
{\rm s.t.}
&~~|{\rm supp}(x)|\le K,
\end{align*}
where $A^T=(a_1,\ldots,a_m)$, $b=(b_1,\ldots,b_m)^T$,
and $K$ is an integer with $1\le K\le n$.
When we convert this problem to problem (P), lower bounds and upper bounds
on $x$, i.e. $L\le x_i \le U$, can be imposed for a
sufficiently large positive number $U$ and a sufficiently small
negative number $L$.

Cardinality constrained linear-quadratic optimal control was investigated in \cite{GaoLi2011}.
Furthermore, a polynomially solvable case of the cardinality-constrained quadratic
optimization problem was identified in \cite{gaoli2013a}.

We focus in this paper on exact solution methods for problem (P).
Standard MIQP solvers that are based on branch-and-bound frameworks can be
applied to (P) directly.
However, the lower bound generated from the continuous relaxation of (P) by relaxing
$y_i\in\{0,1\}$ to $y_i\in [0,1]$ is often quite loose.
Equivalent reformulations with tighter continuous relaxation, i.e., a larger lower bound,
have been proposed
in the literature
\cite{frangioni06mp, frangioni2007letter, frangioni2009letter, zheng2013}.
These reformulations are more efficient when solved in MIQP solvers.
We propose a lift-and-convexification approach to construct a tight reformulation
for problem (P). This lift-and-convexification approach lifts the quadratic term involving $x$ only in the original objective function $f(x,y)$ to a quadratic function
of both $x$ and $y$ and convexifies this equivalent objective function in a quadratic form of $(x,y)$.
The new reformulation retains the linearly constrained structure of the MIQP form so that its
continuous relaxations can be solved efficiently.
At the same time, the lower bound achieved by the continuous relaxation of this newly proposed reformulation
is the same as the lower bound obtained from the state-of-the-art perspective
reformulation.
Thus it improves the performance of branch-and-bound algorithms by providing the same
tightness at the root node as the state-of-the-art perspective reformulation and
much faster processing time at children nodes.
We then further combine our lift-and-convexification approach and the
quadratic convex reformulation (QCR) \cite{billionnet2008, billionnet2009qcr} approach in the literature
to form an even tighter reformulation.
The QCR approach has been applied to
zero-one quadratic programs \cite{billionnet2009qcr} and
integer quadratic programs \cite{billionnet2012qcrExtended}.
While the QCR approach cannot be directly applied to problem (P),
it can be successfully applied on top of our new lift-and-convexification reformulation.
This further reduces the duality gap as we will show in our
numerical tests.

The paper is organized as follows: In \S 2, we review the current state-of-the-art
reformulation and exact solution methods for problem (P).
In \S 3, we propose a lift-and-convexification approach to obtain
a tight reformulation.
We show that this new reformulation is as tight as the state-of-the-art
reformulation in terms of the lower bound from its continuous relaxation.
As the continuous relaxation of the new reformulation is a quadratic program, it
can be thus solved efficiently.
In \S 4, we conduct numerical experiments to demonstrate the
effectiveness of our new lift-and-convexification reformulation.
In \S 5, we review the QCR approach in the literature for
the binary and integer quadratic programs.
We then combine lift-and-convexification approach and the
QCR approach
to form an even tighter reformulation.
We conclude our paper in \S 6.

{\it Notation:} Throughout this paper, we denote by $v(\cdot)$ the optimal value of
problem $(\cdot)$, and $\Re^n_+$ the nonnegative orthant of $\Re^n$.
For any $a\in\Re^n$, we denote by ${\rm diag}(a)={\rm diag}(a_1,\ldots,a_n)$ the diagonal
matrix with $a_i$ being its $i$th diagonal element.
We denote by $e$ the all-one vector.

\section{Literature review and related work}
One efficient solution method for (P) is the
perspective reformulation proposed by \cite{frangioni06mp, frangioni2007letter}, in which problem (P) is transformed into the following equivalent form:
\begin{align*}
({\rm PR}(\rho))~~
\min
&~~f_{\rho}(x,y)=x^{T}(Q-{\rm diag}(\rho))x+c^Tx+h^Ty
+\sum^{n}_{i=1}[\rho_i(x^2_i/y_i)] \\
{\rm s.t.}
&~~Ax+By\le d,\non\\
&~~a_iy_i\le x_i\le b_i y_i,~y_i\in\{0,1\},~i=1,\ldots,n,
\end{align*}
where $\rho\in\Re^{n}$ is chosen such that
\begin{align*}
\rho\geq0~{\rm and}~Q-{\rm diag}(\rho)\succeq0,
\end{align*}
with an assumption $0/0=0$.

The perspective  reformulation is very tight, i.e., the lower bound generated from the
continuous relaxation of this reformulation is usually much higher than
the lower bound generated directly from the continuous relaxation of (P).
To deal with the fractional terms in the objective function of $({\rm PR}(\rho))$,
two tractable reformulations of $({\rm PR}(\rho))$ were proposed in the literature.

The first reformulation is a second-order cone programming (SOCP) reformulation
\cite{akturk2009strong, gunluk2010perspective}.
For each $i$, introducing an additional variable $\phi_i=x_i^2/y_i$ and then rewriting the
constraint $\phi_i\geq x_i^2/y_i$ as an SOCP constraint yields the following SOCP reformulation:
\begin{align*}
({\rm SOCP}(\rho))~~
\min
&~~x^{T}(Q-{\rm diag}(\rho))x+c^Tx+h^Ty+\rho^T\phi \\
{\rm s.t.}
&~~Ax+By\le d,\non\\
&~~a_iy_i\le x_i\le b_i y_i,~y_i\in\{0,1\},~i=1,\ldots,n, \\
&~\left\| \begin{array}{c}
                  x_i \\
                  \frac{\phi_i-y_i}{2} \\
                 \end{array}
   \right\| \le\frac{\phi_i+y_i}{2},~i=1,\ldots,n.
\end{align*}
However, as the problem size grows, the time needed to solve the above SOCP relaxation becomes
a critical factor. When interior point methods are used,
the corresponding branch-and-bound algorithm may converge very slowly.

The second reformulation is the perspective cut (PC) reformulation
\cite{frangioni06mp, frangioni2007letter}.
Representing the value of $x_i^2/y_i$ by the supremum of a set of infinitely many
hyperplanes, which are called perspective cuts, gives rise to the following
PC reformulation:
\begin{align}
({\rm PC}(\rho))~~
\min
&~~x^{T}(Q-{\rm diag}(\rho))x+c^Tx+h^Ty+\rho^T\phi \non \\
{\rm s.t.}
&~~Ax+By\le d,\non\\
&~~a_iy_i\le x_i\le b_i y_i,~y_i\in\{0,1\},~i=1,\ldots,n, \non \\
&~~\phi_i\geq 2\ovl{x}_ix_i-\ovl{x}_i^2y_i,
         \forall \ovl{x}_i\in [a_i,b_i],~i=1,\ldots,n.
         \label{pcut}
\end{align}
The perspective cuts in (\ref{pcut}) can be added dynamically when ${\rm(PC(\rho))}$
is solved in a branch-and-cut framework
(see \cite{frangioni2009letter}).
With the help of warm start and dual methods, quadratic programming relaxations in
the perspective cut algorithm can be solved efficiently.

Let $({\rm\ovl{PR}}(\rho))$, $({\rm\ovl{SOCP}}(\rho))$ and $({\rm\ovl{PC}}(\rho))$
denote the continuous relaxations of\\ $({\rm PR}(\rho))$, $({\rm SOCP}(\rho))$ and
$({\rm PC}(\rho))$, respectively, by relaxing $y_i\in\{0,1\}$ to $y_i\in [0,1]$.
It is easy to see that the objective values of these continuous relaxations form the same
lower bound for $({\rm P})$.
A key issue is how to choose the vector $\rho$ such that this lower bound is as large as possible.
One natural way is to set every component of $\rho$ to be the smallest eigenvalue of $Q$.
Frangioni and Gentile \cite{frangioni2007letter} proposed a better heuristic and set $\rho$ to be
the optimal solution to the following SDP problem:
\begin{align}
\max\{e^T\rho\mid\rho\geq0,Q-{\rm diag}(\rho)\succeq0\}.
\end{align}
Ideally, the best parameter $\rho$ that maximizes the lower bound
$v({\rm\ovl{PR}}(\rho))$ can be found by solving the
following problem:
\begin{align}
({\rm MAX}\rho)~~\max\{v({\rm\ovl{PR}}(\rho))\mid\rho\geq0,Q-{\rm diag}(\rho)\succeq0\}.
\label{maxrho}
\end{align}
Recently, Zheng et al. \cite{zheng2013} established the following interesting result.
\begin{thm}
\label{thmsdpl}
Problem $({\rm MAX}\rho)$ is equivalent to the following
semi-definite programming (SDP) problem:
\begin{align}
({\rm SDP}_l)~~\non
\max
&~~\tau \non\\
{\rm s.t.}
&~ \left(
    \begin{array}{cc}
    \rho_i+\mu_i                                & \frac{1}{2}(c_i-\lambda_i-(a_i+b_i)\mu_i) \\
    \frac{1}{2}(c_i-\lambda_i-(a_i+b_i)\mu_i)^T & h_i-\pi_i+(B^T\eta)_i+\mu a_ibi
    \end{array}
    \right)\succeq 0,\non\\
&~~i=1,\ldots,n,\non\\
&~ \left(
         \begin{array}{cc}
         Q - {\rm diag}(\rho)           & \frac{1}{2}(\lambda+A^T\eta) \\
         \frac{1}{2}(\lambda+A^T\eta)^T & -\eta^Td-e^T\pi-\tau
         \end{array}
    \right)\succeq 0,
    \non\\
&~~(\eta,\mu,\pi,\rho)\in\Re^m_+\times\Re^n_+\times\Re^n_+\times\Re^n_+,\non \\
&~~(\lambda,\tau)\in \Re^n\times \Re. \non
\end{align}
\end{thm}
Zheng et al. \cite{zheng2013} showed that the perspective cut approach for $({\rm PC}(\rho^*))$ with
$\rho^*$ obtained from $({\rm SDP}_l)$ is most efficient  for solving problem (P)
to its optimality.

When $B\equiv 0$ in the constraint $Ax+By\le d$, Frangioni et al. \cite{frangioni2011or} developed an equivalent MIQP reformulation of $({\rm PR}(\rho))$, whose continuous relaxation becomes a quadratic programming problem.
Frangioni et al. \cite{frangioni2013} also proposed an MIQP reformulation of the original problem (P).
But the continuous relaxation of this MIQP reformulation is in general not as tight as
that of the perspective reformulation.

\section{Lift-and-convexification approach}
In this section, we derive a tight \\MIQP reformulation of (P)
by proposing
a lift-and-convexification approach.
This approach  lifts the quadratic term involving $x$ only in the original objective function $f(x,y)$
to an equivalent quadratic function of both $x$ and $y$ and convexifies this equivalent objective function in a quadratic form of $(x,y)$.

Contrast to the the perspective reformulation which involves fractional terms, our new reformulation is a quadratic programming problem whose
continuous relaxations can be solved efficiently.
At the same time, the lower bound achieved by the continuous relaxation of this new reformulation
can be proved to achieve the same lower bound obtained from $({\rm\ovl{SOCP}}(\rho^*))$ or
${\rm(\ovl{PC}(\rho^*))}$ with $\rho^*$ calculated from $({\rm SDP}_l)$.
To construct the new reformulation, we only need to solve an additional SOCP problem, given the
solution for $({\rm SDP}_l)$.

Let us determine first what kind of quadratic functions in the $(x,y)$-space
we need to add to achieve the above mentioned goals.
\smallskip
\begin{thm}
Let $q(x,y)$ be a quadratic function of  $x$ and $y$.
If $q(x,y)=0$ for all $(x,y)\in\{(x,y)\mid
a_iy_i\le x_i\le b_iy_i$, $y_i\in\{0,1\}$,~$i=1,\ldots,n\}$,
where $a_i<b_i,~i=1,~\ldots, ~n$, then $q(x,y)$ must take
the following form:
\begin{align}
q(x,y)=\sum_{i=1}^nq_i(x_i,y_i),
\label{q_1}
\end{align}
where
\begin{align}
q_i(x_i,y_i)=u_i x_iy_i+v_i y_i^2-u_ix_i-v_iy_i
\label{q_2}
\end{align}
is a quadratic function of $(x_i,y_i)$ parameterized by $(u_i,v_i).$
\label{thm_qform}
\end{thm}
\bigskip
\proof
If $(x,y)\in\{(x,y)\mid
a_iy_i\le x_i\le b_iy_i$, $y_i\in\{0,1\}$,~$i=1,\ldots,n\}$,
then for any $i=1,..,n$, if $y_i=0$, then $x_i$ must be $0$.

Let $q(x,y)$ be of the following general form:
\begin{align*}
  q(x,y)
  = (x^T,y^T)
      \left( \begin{array}{cc}
             P   & F \\
             F^T & G \\
             \end{array}
      \right)
      \left( \begin{array}{c}
             x \\
             y \\
             \end{array}
      \right)
      +p^Tx+g^Ty,
\end{align*}
parameterized by $(P,F,G,p,g)\in\mathbb{S}^n\times\Re^{n\times n}\times\mathbb{S}^n
\times\Re^n\times\Re^n$.
Let $x^i$ be a vector with a non-zero component only in its $i$th position and
$y^i$ be a vector with a non-zero component (which is set at one) only in its $i$th position.
It is clear that $(x^i,y^i)\in\{(x,y)\mid
a_iy_i\le x_i\le b_iy_i$, $y_i\in\{0,1\}$,~$i=1,\ldots,n\}$.
If $q(x,y)=0$ for all $(x,y)\in\{(x,y)\mid
a_iy_i\le x_i\le b_iy_i$, $y_i\in\{0,1\}$,~$i=1,\ldots,n\},$
 then
\begin{align}
q(x^i,y^i)        &=0,\label{q_3}\\
q(x^i+x^j,y^i+y^j)&=0,\label{q_4}
\end{align}
for any $i,j=1,..,n$.
(\ref{q_3}) implies
\begin{align}
P_{ii}x_i^2+2F_{ii}x_i+G_{ii}+p_ix_i+g_i=0.\label{q41}
\end{align}
Because (\ref{q41}) must hold for any $x_i\in[a_i, b_i]$ and $a_i<b_i$, we have
\begin{align}
P_{ii} &=0,   \label{q_5}\\
2F_{ii}&=-p_i,\label{q_6}\\
G_{ii} &=-g_i.\label{q_7}
\end{align}
Furthermore, (\ref{q_4}) leads to
\begin{align}
P_{ii}x_i^2+2P_{ij}x_ix_j+P_{jj}x_j^2+2((F_{ii}+F_{ij})x_i+(F_{jj}+F_{ji})x_j)\non\\
+G_{ii}+2G_{ij}+G_{jj}+p_ix_i+p_jx_j+g_i+g_j=0,\non
\end{align}
which can be simplified to the following equality by using (\ref{q_5})-(\ref{q_7}),
\begin{align}
P_{ij}x_ix_j+F_{ij}x_i+F_{ji}x_j+G_{ij}=0.\non
\end{align}
As the above equality holds for any $(x_i,x_j)\in[a_i,b_i]\times[a_j,b_j]$
with $a_i<b_i$ and $a_j<b_j$, we must have
$0=P_{ij}=F_{ij}=F_{ji}=G_{ij}.$
Combining the above equality with (\ref{q_5})-(\ref{q_7}) yields
the following form of $q(x,y)$,
\begin{align}
q(x,y)=\sum_{i=1}^n2F_{ii}x_iy_i+G_{ii}y_i^2-2F_{ii}x_i-G_{ii}y_i,\non
\end{align}
which is of the same form as (\ref{q_2}).
\endproof

We propose the following reformulation of (P):
\begin{align*}
({\rm P}(u,v))~~
\min
&~~f_{u,v}(x,y)=f(x,y)+\sum_{i=1}^nq_i(x_i,y_i)\nonumber\\
s.t.
&~~Ax+By\le d,\non\\
&~~a_iy_i\le x_i\le b_i y_i,~y_i\in\{0,1\},~i=1,\ldots,n,
\end{align*}
where $q_i(x_i,y_i)$ is defined in (\ref{q_2}).
It is easy to see that problem $({\rm P}(u,v))$  is equivalent to (P)
and the continuous relaxation of $({\rm P}(u,v))$ is a quadratic program.

The difference among equivalent formulations (P), $({\rm PR}(\rho))$, and $({\rm P}(u,v))$
lies in their objective functions.
The following example shows the relative relationship among these three objective
functions.
\begin{exam}{\rm
Consider a univariate function $f(x)=x^2-4x$, where
$(x,y)\in\Omega=\{(x,y)\mid y\le x \le 3y,~y\in [0,1]\}$.
Let $q(x,y)=-xy+y^2+x-y$.
Then, $q(x,y)$ is zero at the region
$\{(x,y)\mid y\le x \le 3y,~y\in \{0,1\}\}$.
Figure \ref{pc_lift} illustrates the original function $f(x)$, the lifted quadratic
function $f_{u,v}(x,y):=f(x)+q(x,y)$ and the perspective function $f_p(x,y)=\frac{x^2}{y}-4x$.
While the three functions has the same value when $y\in\{0,1\}$,
we can see that $f_{u,v}(x,y)$  always lies between
$f(x)$ and $f_p(x,y)$ in the region $\Omega$.
This can be numerically verified since
$f_{u,v}(x,y)-f(x)=(x-y)(1-y)\ge0$ for all $(x,y)\in\{(x,y)\mid y\le x \le 3y,~y\in [0,1]\}$.
We can also verify that $f_p(x,y)-f_{u,v}(x,y)=(1-y)y(x^2+y^2-xy)\ge0$
for all $(x,y)\in\{(x,y)\mid y\le x \le 3y,~y\in (0,1]\}$.
\begin{figure}[ht!]
\begin{center}
\includegraphics[width=4in]{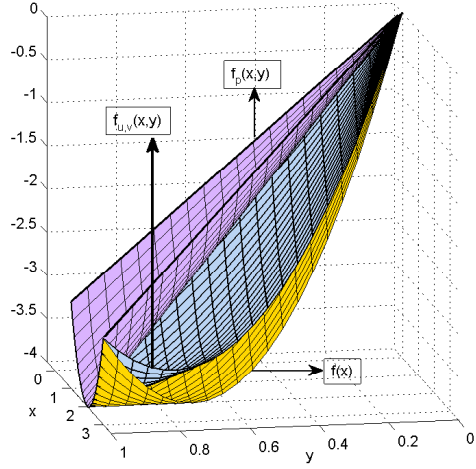}
\caption{\label{pc_lift}  Illustration of the lifting function}
\end{center}
\end{figure}
}
\label{exam1}
\end{exam}

Example \ref{exam1} shows that the objective function of our
new reformulation could lie below that of the perspective reformulation.
This indicates that our new reformulation may  not attain
a lower bound that is tighter than the perspective reformulation.
However, we will show that our new reformulation can achieve the same
lower bound as the perspective reformulation.

Let $u=(u_1,...,u_n)^T$ and $v=(v_1,...,v_n)^T$.
Now, one critical question is ``What is the best parameter vector of $(u,v)$?''
Let $({\rm\ovl{P}}(u,v))$ denote the continuous relaxation of $({\rm P}(u,v))$ by relaxing
$y_i\in\{0,1\}$ to $y_i\in [0,1]$.
It is desirable to choose $(u,v)$ such that the continuous relaxation of $({\rm P}(u,v))$
is as tight as possible.
This clear goal motivates us to consider the following problem:
\begin{align}
({\rm MAX}uv)~~\max\{v({\rm\ovl{P}}(u,v))\mid u,v\in\Re^n,~f_{u,v}(x,y)~is~convex\}.
\label{maxq}
\end{align}

\begin{thm}
\label{thmsdpq}
Problem $({\rm MAX}uv)$ is equivalent to the following SDP problem:
\begin{align}
({\rm SDP}_q)~~
\max
&~~\tau \non \\
{\rm s.t.}
&~  \left(
    \begin{array}{ccc}
    Q                                      & \frac{1}{2}{\rm diag}(u)                          & \frac{1}{2}\alpha(u,\eta,\mu,\sigma) \\
    \frac{1}{2}{\rm diag}(u)               & {\rm diag}(v)                                     & \frac{1}{2}\beta(v,\eta,\mu,\sigma,\lambda,\pi)\\
    \frac{1}{2}\alpha(u,\eta,\mu,\sigma)^T & \frac{1}{2}\beta(v,\eta,\mu,\sigma,\lambda,\pi)^T & -\eta^T d-e^T\pi-\tau \\
    \end{array}
    \right)\succeq 0,
    \label{sdpQ1}\\
&~~(\eta,\mu,\sigma,\lambda,\pi)\in\Re^m_+\times\Re^n_+\times\Re^n_+\times\Re^n_+\times\Re^n_+,
   \label{sdpQ2} \\
&~~(u,v,\tau)\in \Re^n\times \Re^n\times \Re, \non
\end{align}
where
\begin{align}
  \alpha(u,\eta,\mu,\sigma)&=c-u+A^T\eta-\mu+\sigma,\label{dsc1}\\
  \beta(v,\eta,\mu,\sigma,\lambda,\pi)&=h-v+B^T\eta+{\rm diag}(a)\mu-{\rm diag}(b)\sigma-\lambda+\pi.\label{dsc2}
\end{align}
\end{thm}

\proof
We first express $({\rm\ovl{P}}(u,v))$ by its dual form.
Associate the following multipliers to the constraints in $({\rm\ovl{P}}(u,v))$:
\begin{itemize}
  \item $\eta\in\Re^m_+$ for $Ax+By\le d$;
  \item $\mu_i$ and $\sigma_i\in \Re_+$ for
        $a_iy_i\le x_i$ and $x_i\le b_iy_i$, respectively, $i=1,\ldots,n$; and
  \item $\lambda_i$ and $\pi_i\in \Re_+$ for $y_i\ge 0$ and $y_i\le 1$, respectively,
        $i=1,\ldots,n$.
\end{itemize}
Let $\mu=(\mu_1,\ldots,\mu_n)^T$, $\sigma=(\sigma_1,\ldots,\sigma_n)^T$,
$\lambda=(\lambda_1,\ldots,\lambda_n)^T$ and $\pi=(\pi_1,\ldots,\pi_n)^T$.
Let $\omega$ denote the vector formed by $\eta,\mu,\sigma,\lambda$ and $\pi$.
The Lagrangian function of $({\rm\ovl{P}}(u,v))$ is then given by
\begin{align*}
  L(x,y,\omega)
  &= x^TQx+c^Tx+h^Ty+\sum_{i=1}^n(u_ix_iy_i+v_iy_i^2-u_ix_i-v_iy_i)\\
  &~~~+\eta^T(Ax+By-d)+\sum_{i=1}^{n}\mu_i(a_iy_i-x_i)+\sum_{i=1}^{n}\sigma_i(x_i-b_iy_i)\\
  &~~~+\sum_{i=1}^{n}\lambda_i(-y_i)+\sum_{i=1}^{n}\pi_i(y_i-1)\\
  &= (x^T,y^T)
      \left( \begin{array}{cc}
             Q                        & \frac{1}{2}{\rm diag}(u) \\
             \frac{1}{2}{\rm diag}(u) & {\rm diag}(v) \\
             \end{array}
      \right)
      \left( \begin{array}{c}
             x \\
             y \\
             \end{array}
      \right)
      +(c-u+A^T\eta-\mu+\sigma)^{T}x\\
  &~~~+(h-v+B^T\eta+{\rm diag}(a)\mu-{\rm diag}(b)\sigma-\lambda+\pi)^Ty
      -\eta^T d-e^{T}\pi.
\end{align*}
Furthermore, the Lagrangian dual problem of $({\rm\ovl{P}}(u,v))$ can be expressed as
\begin{align}
\max \{\min_{(x,y)\in\Re^n\times \Re^n} ~L(x,y,\omega)\mid \omega\ge 0\}.\label{ds1}
\end{align}
Introducing an additional variable $\tau$, we can rewrite (\ref{ds1}) as
\begin{align}
\max
&~~\tau\label{ds21}\\
{\rm s.t.}
&~~\min_{(x,y)\in\Re^n\times \Re^n} ~L(x,y,\omega)\ge\tau,\label{ds22}\\
&~~\omega\ge 0.\label{ds23}
\end{align}
We see that the constraint in (\ref{ds22}) is equivalent to
$L(x,y,\omega)-\tau\ge 0$ for all $x,y$, which is further equivalent to
\begin{align}
  L(x/t,y/t,\omega)-\tau\ge 0,
  ~~\forall (x,y)\in\Re^n\times \Re^n,~\forall t\in\Re,~t\not=0.\label{Lhomo1}
\end{align}
Multiplying both sides of (\ref{Lhomo1}) by $t^2$ yields a homogeneous quadratic form
of $(x,y,t)$ in the left-hand side of (\ref{Lhomo1}).
Thus, the constraint in (\ref{ds22}) is equivalent to the following semidefinite constraint:
\begin{align}
  \left(
  \begin{array}{ccc}
  Q                                      & \frac{1}{2}{\rm diag}(u)                          & \frac{1}{2}\alpha(u,\eta,\mu,\sigma) \\
  \frac{1}{2}{\rm diag}(u)               & {\rm diag}(v)                                     & \frac{1}{2}\beta(v,\eta,\mu,\sigma,\lambda,\pi)\\
  \frac{1}{2}\alpha(u,\eta,\mu,\sigma)^T & \frac{1}{2}\beta(v,\eta,\mu,\sigma,\lambda,\pi)^T & -\eta^T d-e^T\pi-\tau \\
  \end{array}
  \right)\succeq 0,\label{ds31}
\end{align}
where $\alpha(u,\eta,\mu,\sigma)$ and $\beta(v,\eta,\mu,\sigma,\lambda,\pi)$ are defined
by (\ref{dsc1}) and (\ref{dsc2}), respectively.
Consequently, the problem in (\ref{ds21})-(\ref{ds23}) can be expressed as
\begin{align}
\max\{\tau\mid {\rm constraint}~(\ref{ds31}),~\omega\ge 0\}.\label{ds4}
\end{align}
If $f_{u,v}(x,y)$ is convex, by the strong duality of convex quadratic programming
(see, e.g., Proposition 6.5.6 in \cite{bertsekas2003}),
the optimal values of $({\rm\ovl{P}}(u,v))$ and (\ref{ds4}) are equal.
Thus,  we have shown that problem $({\rm MAX}uv)$ is equivalent to an SDP problem
in the form of $({\rm SDP}_q)$.
\endproof

Let $(u^*,v^*)$ be the optimal parameters for our new reformulation $({\rm P}(u,v))$
by solving the SDP program $({\rm SDP}_q)$.
Let $\rho^*$ be the optimal parameters for $({\rm PR}(\rho))$
by solving the SDP program $({\rm SDP}_l)$.
It is necessary to compare the tightness of $v({\rm\ovl{P}}(u^*,v^*))$ and $v({\rm\ovl{PR}}(\rho^*))$,
i.e., the bounds from the continuous relaxation of the ``best'' reformulation proposed in this paper and the ``best'' perspective reformulation.
We will show in the following that these two bounds are the same  using constructive proofs.
\smallskip
\begin{thm}
Define $\ovl{\rho}\in\Re^n$ as
\begin{align}
   \ovl{\rho}_i &=  \begin{cases}
                    0                       &{\rm~if~}v^*_i=0,  \\
                    \frac{u^{*2}_i}{4v^*_i} &{\rm~otherwise},  \\
                    \end{cases}~~i=1,\ldots,n.
\label{findrho}
\end{align}
Then,
\begin{description}
\item[(a)] $\ovl{\rho}$ is feasible for problem $({\rm MAX}\rho)$,
\item[(b)] $v({\rm\ovl{P}}(u^*,v^*))\le v({\rm\ovl{PR}}(\ovl{\rho}))$.
\end{description}
\label{thm_rho_bar}
\end{thm}
\bigskip
\proof
(a)
From the proof of Theorem \ref{thmsdpq}, we know that there exists
$(\tau^*,\eta^*,\mu^*,\sigma^*,$\\$\lambda^*,\pi^*)$ such that
$(u^*,v^*,\tau^*,\eta^*,\mu^*,\sigma^*,\lambda^*,\pi^*)$ is optimal to $({\rm SDP}_q)$.
The constraint in (\ref{sdpQ1}) implies that
\begin{align*}
   \left( \begin{array}{cc}
          Q                          & \frac{1}{2}{\rm diag}(u^*) \\
          \frac{1}{2}{\rm diag}(u^*) & {\rm diag}(v^*) \\
          \end{array}
   \right) \succeq 0.
\end{align*}
For any $g\in\Re^n$, we define $p\in\Re^n$ with
\begin{align*}
        p_i   &=  \begin{cases}
                  0                        &{\rm~if~}v^*_i=0,  \\
                  -\frac{u^*_i}{2v^*_i}g_i &{\rm~otherwise},  \\
                  \end{cases}~~i=1,\ldots,n.
\end{align*}
We then have
\begin{align*}
       g^T(Q-\ovl{\rho})g
  &=    g^TQg+\sum^n_{i=1}[g_i^2(-\ovl{\rho}_i)] \\
  &=    g^TQg+\sum^n_{i=1}[u^*_ig_ip_i+v^*_ip_i^2] \\
  &=    g^TQg+g^T{\rm diag}(u^*)p+p^T{\rm diag}(v^*)p \\
  &=     (g^T,p^T)
         \left( \begin{array}{cc}
         Q                          & \frac{1}{2}{\rm diag}(u^*) \\
         \frac{1}{2}{\rm diag}(u^*) & {\rm diag}(v^*) \\
         \end{array}
         \right)
         \left( \begin{array}{c}
                g \\
                p \\
                \end{array}
         \right) \\
  &\geq 0.
\end{align*}
Hence $(Q-\ovl{\rho})\succeq 0$ and $\ovl{\rho}$ is feasible to problem $({\rm MAX}\rho)$.

(b) As $({\rm\ovl{PR}}(\ovl{\rho}))$ and $({\rm\ovl{P}}(u^*,v^*))$
have the same feasible region, for any feasible $(x,y)$,
we can compare their objective values as follows,
\begin{align}
  &f_{\rho}(x,y;~\ovl{\rho}) - f_{u,v}(x,y;~u^*,v^*)
         \label{compare} \\
  &=     \sum^{n}_{i=1}[
         \ovl{\rho}_i(\frac{x^2_i}{y_i})-\ovl{\rho}_ix^2_i
         -(u^*_ix_iy_i+v^*_iy_i^2-u^*_ix_i-v^*_iy_i)
         ] \non \\
  &=     \sum^{n}_{i=1}[
         \frac{1}{4v^*_i}
         (u_i^{*2}x_i^2\frac{1-y_i}{y_i}+4u_i^*v_i^*x_i(1-y_i)+4v_i^{*2}y_i(1-y_i))
         ] \non \\
  &=     \sum^{n}_{i=1}[
         \frac{y_i(1-y_i)}{4v^*_i}
         ((\frac{u_i^*x_i}{y_i})^2+4u_i^*v_i^*\frac{x_i}{y_i}+4v_i^{*2})
         ] \non \\
  &=     \sum^{n}_{i=1}[
         \frac{y_i(1-y_i)}{4v^*_i}
         (u^*_i\frac{x_i}{y_i}+2v^*_i)^2
         ] \non \\
  &\geq 0. \non
\end{align}
The above deduction is valid because if $v^*_i=0$, then $u^*_i=0$ due to (\ref{sdpQ1}).
Thus $v({\rm\ovl{P}}(u^*,v^*))\le v({\rm\ovl{PR}}(\ovl{\rho}))$.
\endproof

The following corollary is a direct result of Theorem \ref{thm_rho_bar}.
\smallskip
\begin{cor}
$v({\rm\ovl{P}}(u^*,v^*))\le v({\rm\ovl{PR}}(\rho^*))$.
\label{no_better}
\end{cor}
\bigskip

Next we show the other way around.
\smallskip
\begin{thm}
\label{thmfinduv}
Suppose that $(x^*,y^*)$ is optimal to
problem $({\rm\ovl{PR}}(\rho^*))$, the continuous relaxation of $({\rm PR}(\rho^*))$.
Define $\ovl{u},\ovl{v}\in\Re^n$ with
\begin{align}
   \ovl{u}_i &=  \begin{cases}
                 0                             &{\rm~if~}y^*_i=0,  \\
                 -2\rho^*_i\frac{x^*_i}{y^*_i} &{\rm~otherwise},  \\
                 \end{cases}~~i=1,\ldots,n,
   \label{findu} \\
   \ovl{v}_i &=  \begin{cases}
                 0                                 &{\rm~if~}y^*_i=0,  \\
                 \rho^*_i\frac{x^{*2}_i}{y^{*2}_i} &{\rm~otherwise},  \\
                 \end{cases}~~i=1,\ldots,n.
   \label{findv}
\end{align}
Then,
\begin{description}
\item[(a)] $(\ovl{u},\ovl{v})$ is feasible to problem $({\rm MAX}uv)$,
\item[(b)] $v({\rm\ovl{P}}(\ovl{u},\ovl{v}))=v({\rm\ovl{PR}}(\rho^*))$.
\end{description}
\label{thm_uv_bar}
\end{thm}
\bigskip
\proof
(a)
For any $g,p\in\Re^n$, we have
\begin{align*}
  &      (g^T,p^T)
         \left( \begin{array}{cc}
         Q                              & \frac{1}{2}{\rm diag}(\ovl{u}) \\
         \frac{1}{2}{\rm diag}(\ovl{u}) & {\rm diag}(\ovl{v}) \\
         \end{array}
         \right)
         \left( \begin{array}{c}
                g \\
                p \\
                \end{array}
         \right) \\
  &=    g^TQg+g^T{\rm diag}(\ovl{u})p+p^T{\rm diag}(\ovl{v})p \\
  &=    g^TQg+\sum^n_{i=1}[\ovl{u}_ig_ip_i+\ovl{v}_ip_i^2] \\
  &=    g^T(Q-\rho^*)g+\sum^n_{i=1}[\rho^*_ig_i^2-2\rho^*_i\frac{x^*_i}{y^*_i}g_ip_i
         +\rho^*_i\frac{x^{*2}_i}{y^{*2}_i}p_i^2] \\
  &=    g^T(Q-\rho^*)g+\sum^n_{i=1}[\rho^*_i(g_i-\frac{x^*_i}{y^*_i}p_i)^2] \\
  &\geq 0.
\end{align*}
Hence $(\ovl{u},\ovl{v})$ is feasible to problem $({\rm MAX}uv)$.

(b)
We first show that $(x^*,y^*)$ is also an optimal solution for $({\rm\ovl{P}}(\ovl{u},\ovl{v}))$,
we compare the gradients of $f_{\rho}(x,y;~\rho^*)$ and $f_{u,v}(x,y;~\ovl{u},\ovl{v})$
at the point $(x^*,y^*)$. For $i=1,..,n,$ we have
\begin{eqnarray*}
  & & (\nabla f_{u,v}(x^*,y^*;~\ovl{u},\ovl{v}))_i \\
  &=& 2Q_i^Tx^*+c_i+\ovl{u}_iy^*_i-\ovl{u}_i \\
  &=& 2Q_i^Tx^*+c_i-2\rho^*_i\frac{x^*_i}{y^*_i}*y^*_i+2\rho^*_i\frac{x^*_i}{y^*_i} \\
  &=& 2Q_i^Tx^*+c_i-2\rho^*_ix^*_i+2\rho^*_i\frac{x^*_i}{y^*_i} \\
  &=& (\nabla f_{\rho}(x^*,y^*;~\rho^*))_i
\end{eqnarray*}
and
\begin{eqnarray*}
  & & (\nabla f_{u,v}(x^*,y^*;~\ovl{u},\ovl{v}))_{n+i} \\
  &=& h_i+\ovl{u}_ix^*_i+2\ovl{v}_iy^*_i-\ovl{v}_i \\
  &=& h_i-2\rho^*_i\frac{x^*_i}{y^*_i}*x^*_i+2\rho^*_i\frac{x^{*2}_i}{y^{*2}_i}y^*_i
      -\rho^*_i\frac{x^{*2}_i}{y^{*2}_i} \\
  &=& h_i-\rho^*_i\frac{x^{*2}_i}{y^{*2}_i} \\
  &=& (\nabla f_{\rho}(x^*,y^*;~\rho^*))_{n+i}.
\end{eqnarray*}
So $\nabla f_{u,v}(x^*,y^*;~\ovl{u},\ovl{v})=\nabla f_{\rho}(x^*,y^*;~\rho^*)$.
As $(x^*,y^*)$ is assumed to be optimal to problem $({\rm\ovl{PR}}(\rho^*))$,
the directional derivative of $f_{\rho}(x,y;~\rho^*)$  at $(x^*,y^*)$ along any
feasible direction should be non-negative.
Since the feasible regions of ${\rm\ovl{P}}(\ovl{u},\ovl{v})$ and ${\rm\ovl{PR}}(\rho^*)$
are the same, the directional derivative of $f_{u,v}(x,y;~\ovl{u},\ovl{v})$
at $(x^*,y^*)$ along any feasible direction is also non-negative.
So $(x^*,y^*)$ must also be optimal for $({\rm\ovl{P}}(\ovl{u},\ovl{v}))$ because of
the convexity of $f_{u,v}(x,y;~\ovl{u},\ovl{v})$.
(See e.g., Chapter 2.1 of \cite{borwein2006}.)

Finally, we show $v({\rm\ovl{P}}(\ovl{u},\ovl{v}))=v({\rm\ovl{PR}}(\rho^*))$.
Similar to (\ref{compare}), we have
\begin{align*}
  &  f_{\rho}(x^*,y^*;~\rho^*) - f_{u,v}(x^*,y^*;~\ovl{u},\ovl{v}) \\
  &=     \sum^{n}_{i=1}[
         \rho^*_i(\frac{x^{*2}_i}{y^*_i})-\rho^*_ix^{*2}_i
         -(\ovl{u}_ix^*_iy^*_i+\ovl{v}_iy_i^{*2}-\ovl{u}_ix^*_i-\ovl{v}_iy^*_i)
         ] \non \\
  &=     \sum^{n}_{i=1}[
         \rho^*_i(\frac{x^{*2}_i}{y^*_i})-\rho^*_ix^{*2}_i
         -(-2\rho^*_i\frac{x^*_i}{y^*_i}x^*_iy^*_i+\rho^*_i\frac{x^{*2}_i}{y^{*2}_i}y_i^{*2}+2\rho^*_i\frac{x^*_i}{y^*_i}x^*_i-\rho^*_i\frac{x^{*2}_i}{y^{*2}_i}y^*_i)
         ] \non \\
  &=     \sum^{n}_{i=1}[
         \rho^*_i(\frac{x^{*2}_i}{y^*_i})-\rho^*_ix^{*2}_i
         -(-2\rho^*_ix^{*2}_i+\rho^*_ix_i^{*2}+2\rho^*_i\frac{x^{*2}_i}{y^*_i}-\rho^*_i\frac{x^{*2}_i}{y^{*}_i})
         ] \non \\
  &=  0.
\end{align*}
This completes the proof.
\endproof

The following corollary is a direct result of Theorem \ref{thm_uv_bar}.
\smallskip
\begin{cor}
$v({\rm\ovl{P}}(u^*,v^*))\ge v({\rm\ovl{PR}}(\rho^*))$.
\label{no_worse}
\end{cor}

\bigskip
Combining Corollaries \ref{no_better} and \ref{no_worse} yields
the following result.
\smallskip
\begin{thm}
$v({\rm\ovl{P}}(u^*,v^*))=v({\rm\ovl{PR}}(\rho^*))$.
\label{equal}
\end{thm}

\bigskip
Thus the bound from our new reformulation is as good as the bound
from the perspective reformulation.
However, to find $(u^*,v^*)$, we need to solve $({\rm SDP}_q)$
which is an SDP program that has a much larger size than
$({\rm SDP}_l)$.
Our numerical tests show that $({\rm SDP}_q)$ could consume
ten times of the computation time of $({\rm SDP}_l)$.
Fortunately, based on Theorems \ref{thm_rho_bar} and \ref{thm_uv_bar}, the
following corollary becomes evident which reveals the nonnecessity in using $({\rm SDP}_q)$ in the calculation.
\smallskip
\begin{cor}
\begin{description}
\item[(a)] Define $\ovl{\rho}$ as in (\ref{findrho}).
           Then $\ovl{\rho}$ is optimal for problem $({\rm MAX}\rho)$,
\item[(b)] Define $(\ovl{u},\ovl{v})$ as in (\ref{findu}) and (\ref{findv}).
           Then $(\ovl{u},\ovl{v})$ is optimal to problem $({\rm MAX}uv)$.
\end{description}
\label{no_large_sdp}
\end{cor}

\bigskip
Thus, in order to construct our new reformulation $({\rm P}(u^*,v^*))$, we only need to solve
$({\rm SDP}_l)$ first to get the optimal solution $\rho^*$ for $({\rm MAX}\rho)$
and then
solve the SOCP problem $({\rm\ovl{SOCP}}(\rho^*))$ and construct $(u^*,v^*)$
according to (\ref{findu}) and (\ref{findv}).
Note that $({\rm\ovl{PR}}(\rho^*))$ is equivalent to $({\rm\ovl{SOCP}}(\rho^*))$.

\section{Computational results}

In this section, we conduct computational experiments to compare the performance
of the perspective cut reformulation $({\rm PC}(\rho))$ and our new reformulation
$({\rm P}(u,v))$.
To be specific, we compare the performance of standard MIQP solvers between solving the  following two reformulations of problem $({\rm P})$:
\begin{itemize}
   \item $({\rm PC})$: the perspective reformulation $({\rm PC}(\rho))$ with $\rho=\rho^*$,
   where $\rho^*$ is computed by solving $({\rm SDP}_l)$.
   \item $({\rm LCR})$: our lift-and-convexification reformulation $({\rm P}(u,v))$ with $(u,v)=(u^*,v^*)$,
   where $(u^*,v^*)$ is obtained by first solving $({\rm\ovl{SOCP}}(\rho^*))$ and
   then configuring $(u^*,v^*)$ according to (\ref{findu}) and (\ref{findv}) in
   Theorem \ref{thmfinduv}.
\end{itemize}

Although the continuous relaxation of (LCR) is as tight as that of (PC)
at the root node of the branch-and-bound tree, the relaxations in
(LCR)  are in general looser than those in (PC) at children nodes.
The advantage of (LCR) is that its continuous relaxations are
quadratic programs and thus can be solved much
faster than the continuous relaxations of (PC).
We would like to test if this advantage of (LCR) would dominate (at least verifying itself as a
competitive and useful reformulation).

The time difference between finding $\rho^*$ and $(u^*,v^*)$ is
the time needed to solve one SOCP programming problem $({\rm\ovl{SOCP}}(\rho^*))$.
We will count this amount of time into the computational time for (LCR) in the comparison,
although this amount of time is quite small in general.

The two reformulations and $({\rm\ovl{SOCP}}(\rho^*))$ are all solved in 64-bit
IBM ILOG CPLEX Optimization Studio 12.3
(Hereinafter referred to as \texttt{CPLEX}) through its \texttt{C} interface.
The perspective cut reformulation $({\rm PC})$ is implemented by means of
\texttt{user cut callbacks} and \texttt{lazy constraint callbacks} in \texttt{CPLEX}.
Although \cite{frangioni06mp} suggested to apply the separation procedure
only once at each node, we do not limit the times of separation because we find
that in our numerical tests, if we allow \texttt{CPLEX} to actively generate
mixed integer cuts, the computation would be much faster if the times of separation are unlimited at each node.
$({\rm SDP}_l)$ is solved using \texttt{sedumi} interfaced
by \texttt{CVX} 1.21 (\cite{cvx, gb08}) on \texttt{Matlab R2012b}.

All the computation is conducted on a Linux machine (64-bit CentOS Release 5.5)
with 48 GB of RAM. All the tests are confined on one single thread (2.99 GHz).

We consider two types of test problems from the cardinality constrained mean-variance portfolio selection (MV)
and the subset selection problem (SSP) in our computational experiments.

\subsection{Cardinality constrained portfolio selection problem}

In this subsection, we compare (PC) and (LCR) for the cardinality constrained mean-variance portfolio selection problem
(MV) introduced in the introduction section.

Frangioni and Gentile \cite{frangioni2007letter} tested 90 instances of (MV) in their paper,
$30$ instances each for $n=200$, $300$ and $400$.
The $30$ instances for each $n$ are divided further into three subsets denoted by
$n^-$, $n^0$ and $n^+$, $10$ in each subset, with different diagonal dominance in the
matrix $Q$. We use these 90 instances
created in \cite{frangioni2007letter} in our test.
While Frangioni and Gentile \cite{frangioni2007letter} did not consider the cardinality constraint in their models,
we add the cardinality constraint to these instances in our numerical experiments.
Testing each instance without the cardinality constraint and with $K=6$, $8$, $10$, and $12$,
we have 450 instances of (MV).
The data files of these instances can be downloaded at:
\texttt{http://www.di.unipi.it/optimize/Data/MV.html}.

\begin{table}[t!] \footnotesize \vspace*{0.0in} \centering
\caption{\label{tab1}Numerical results of reformulations for {\rm (MV)}}
\tabcolsep=8pt
\begin{tabular}{cccccccccccc}
\hline \noalign{\smallskip}
     \multirowcell{3}{(MV)}
   & \multirowcell{3}{$K$}
   & \multirowcell{3}{$({\rm time}_l)$}
   & \multirowcell{3}{$({\rm time}_s)$}
   &&\multicolumn{2}{c}{${\rm (PC)}$}
   &&\multicolumn{2}{c}{${\rm (LCR)}$}\\
     \noalign{\smallskip}
     \cline{6-7} \cline{9-10}
     \noalign{\smallskip}
   &&&
   && time & nodes
   && time & nodes \\
\hline \noalign{\smallskip}
   \multirowcell{5}{$200^+$}
   & 6 & 27.92  & 2.30  &  & 19.68  & 65  &  & 4.10  & 26  \\
   & 8 & 27.02  & 2.13  &  & 11.20  & 55  &  & 2.29  & 19  \\
   & 10 & 26.19  & 2.06  &  & 6.18  & 50  &  & 2.86  & 42  \\
   & 12 & 26.29  & 2.05  &  & 10.40  & 111  &  & 4.28  & 95  \\
   & nonK & 32.60  & 2.17  &  & 13.31  & 148  &  & 8.66  & 147  \\
\hline \noalign{\smallskip}
   \multirowcell{5}{$200^0$}
   & 6 & 26.66  & 2.54  &  & 18.03  & 86  &  & 7.07  & 65  \\
   & 8 & 27.20  & 2.05  &  & 16.02  & 96  &  & 7.59  & 73  \\
   & 10 & 23.83  & 2.21  &  & 9.38  & 126  &  & 4.16  & 120  \\
   & 12 & 24.88  & 2.10  &  & 30.73  & 256  &  & 24.44  & 217  \\
   & nonK & 30.64  & 2.00  &  & 33.76  & 291  &  & 27.42  & 281  \\
\hline \noalign{\smallskip}
   \multirowcell{5}{$200^-$}
   & 6 & 26.00  & 2.39  &  & 26.16  & 204  &  & 18.49  & 248  \\
   & 8 & 25.53  & 2.16  &  & 22.49  & 235  &  & 19.96  & 287  \\
   & 10 & 24.67  & 2.03  &  & 19.88  & 328  &  & 11.00  & 350  \\
   & 12 & 25.23  & 2.11  &  & 287.60  & 3306  &  & 329.02  & 1964  \\
   & nonK & 28.47  & 1.85  &  & 152.45  & 1935  &  & 219.22  & 1380  \\
\hline \noalign{\smallskip}
   \multirowcell{5}{$300^+$}
   & 6 & 62.05  & 5.91  &  & 82.11  & 127  &  & 14.25  & 26  \\
   & 8 & 65.26  & 5.54  &  & 48.48  & 133  &  & 8.43  & 25  \\
   & 10 & 55.55  & 5.47  &  & 16.85  & 76  &  & 4.20  & 30  \\
   & 12 & 60.31  & 6.89  &  & 18.43  & 128  &  & 10.27  & 119  \\
   & nonK & 87.59  & 16.23  &  & 108.44  & 446  &  & 26.60  & 241  \\
\hline \noalign{\smallskip}
   \multirowcell{5}{$300^0$}
   & 6 & 60.45  & 5.43  &  & 43.53  & 118  &  & 28.42  & 105  \\
   & 8 & 54.63  & 5.14  &  & 51.06  & 190  &  & 32.20  & 123  \\
   & 10 & 57.15  & 5.04  &  & 22.00  & 148  &  & 11.96  & 129  \\
   & 12 & 61.08  & 5.14  &  & 75.97  & 238  &  & 80.78  & 249  \\
   & nonK & 63.26  & 4.80  &  & 101.63  & 371  &  & 95.79  & 323  \\
\hline \noalign{\smallskip}
   \multirowcell{5}{$300^-$}
   & 6 & 63.27  & 5.10  &  & 55.05  & 236  &  & 48.85  & 237  \\
   & 8 & 62.77  & 5.72  &  & 99.03  & 399  &  & 62.47  & 328  \\
   & 10 & 60.52  & 6.40  &  & 47.20  & 471  &  & 40.47  & 609  \\
   & 12 & 62.08  & 4.92  &  & 35.73  & 312  &  & 48.94  & 341  \\
   & nonK & 68.22  & 4.56  &  & 137.15  & 493  &  & 117.96  & 506  \\
\hline \noalign{\smallskip}
   \multirowcell{5}{$400^+$}
   & 6 & 126.97  & 22.32  &  & 70.89  & 186  &  & 22.73  & 39  \\
   & 8 & 122.96  & 14.19  &  & 196.76  & 587  &  & 22.25  & 38  \\
   & 10 & 122.70  & 12.02  &  & 28.43  & 95  &  & 14.23  & 47  \\
   & 12 & 139.76  & 27.78  &  & 38.70  & 181  &  & 14.37  & 105  \\
   & nonK & 100.40  & 25.47  &  & 562.95  & 849  &  & 364.24  & 613  \\
\hline \noalign{\smallskip}
   \multirowcell{5}{$400^0$}
   & 6 & 104.11  & 11.74  &  & 105.06  & 236  &  & 88.21  & 197  \\
   & 8 & 124.41  & 10.54  &  & 149.88  & 435  &  & 82.51  & 170  \\
   & 10 & 119.51  & 10.47  &  & 54.55  & 276  &  & 49.06  & 287  \\
   & 12 & 128.22  & 10.73  &  & 71.52  & 376  &  & 40.67  & 287  \\
   & nonK & 125.80  & 13.20  &  & 542.91  & 1132  &  & 227.31  & 846  \\
\hline \noalign{\smallskip}
   \multirowcell{5}{$400^-$}
   & 6 & 104.30  & 11.77  &  & 115.01  & 393  &  & 216.02  & 566  \\
   & 8 & 116.33  & 10.08  &  & 341.40  & 1053  &  & 239.26  & 599  \\
   & 10 & 110.30  & 10.08  &  & 83.77  & 515  &  & 95.99  & 564  \\
   & 12 & 118.64  & 10.60  &  & 149.93  & 495  &  & 29.41  & 343  \\
   & nonK & 123.05  & 10.14  &  & 750.27  & 1448  &  & 715.51  & 1373  \\
\hline \noalign{\smallskip}
\end{tabular}
\end{table}

Table \ref{tab1} summarizes the numerical results for the 450 instances of
(MV) when the time limit is set at $10000$ seconds.
Each line reports the average results for the $10$ instances in a subset.
The notations in the table are given as follows: The column ``$({\rm time}_l)$'' is
the computation time for solving $({\rm SDP}_l)$ and the column ``$({\rm time}_s)$''
is the computation time for solving $({\rm\ovl{SOCP}}(\rho^*))$.
The termination threshold of the relative gap
(in percentage) between the objective value of the incumbent solution and the best
lower bound is set to be $0.01\%$.
(The exact value of the relative gap when \texttt{CPLEX} terminates
could range between $0.00\%$ and $0.01\%$. Rounding this number would
make it $0.00\%$ or $0.01\%$).
Because all our instances terminated with a relative gap smaller than
$0.01\%$, the relative gap is not reported here.
The columns ``time'' and ``nodes'' are the computing time (in seconds) and the
number of nodes explored by \texttt{CPLEX} respectively.
The ``nonK'' refers to the instances with no cardinality constraints.

\begin{figure}[ht!]
\begin{center}
\includegraphics[width=5in]{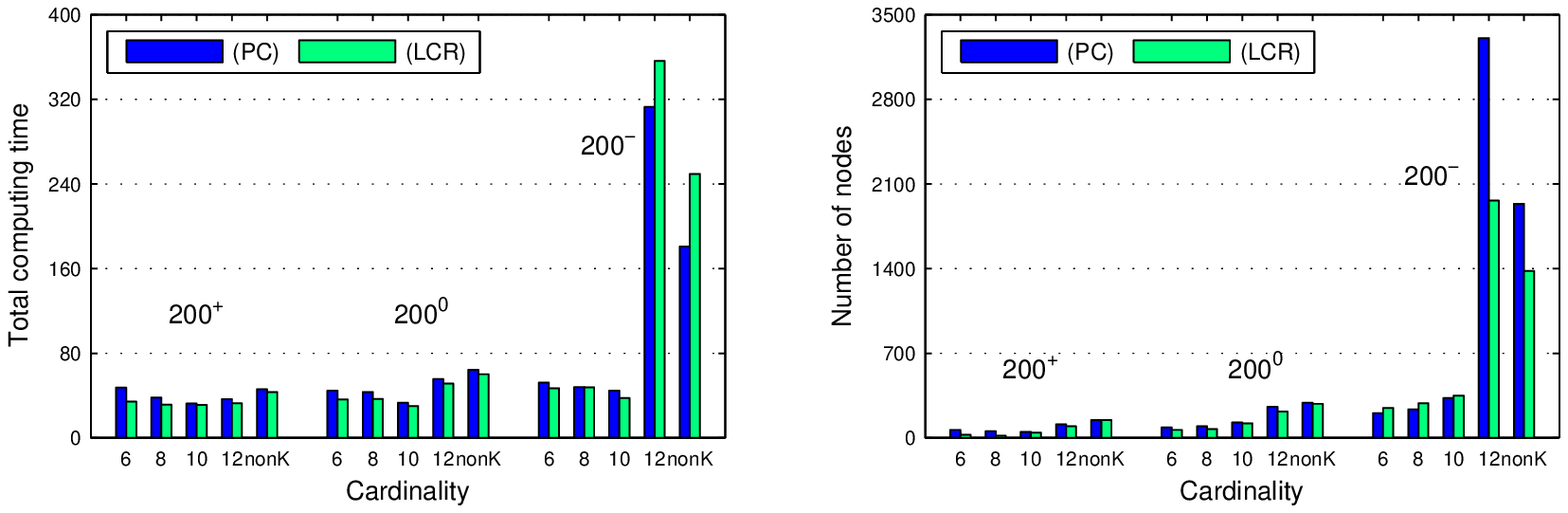}
\includegraphics[width=5in]{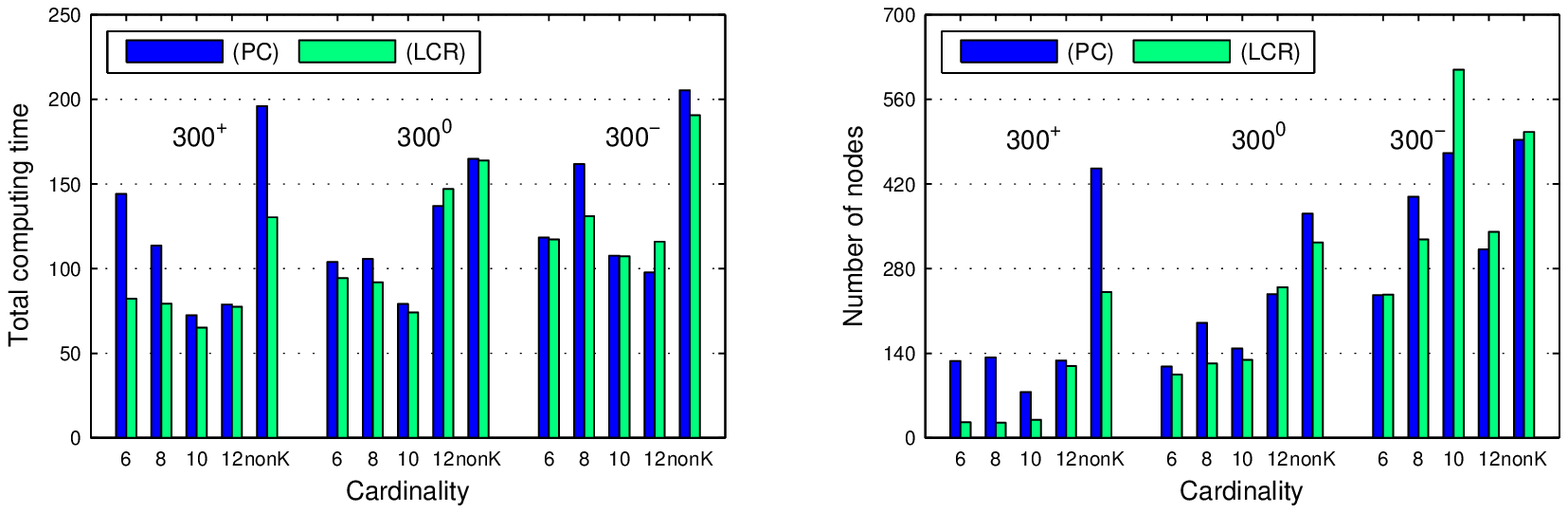}
\includegraphics[width=5in]{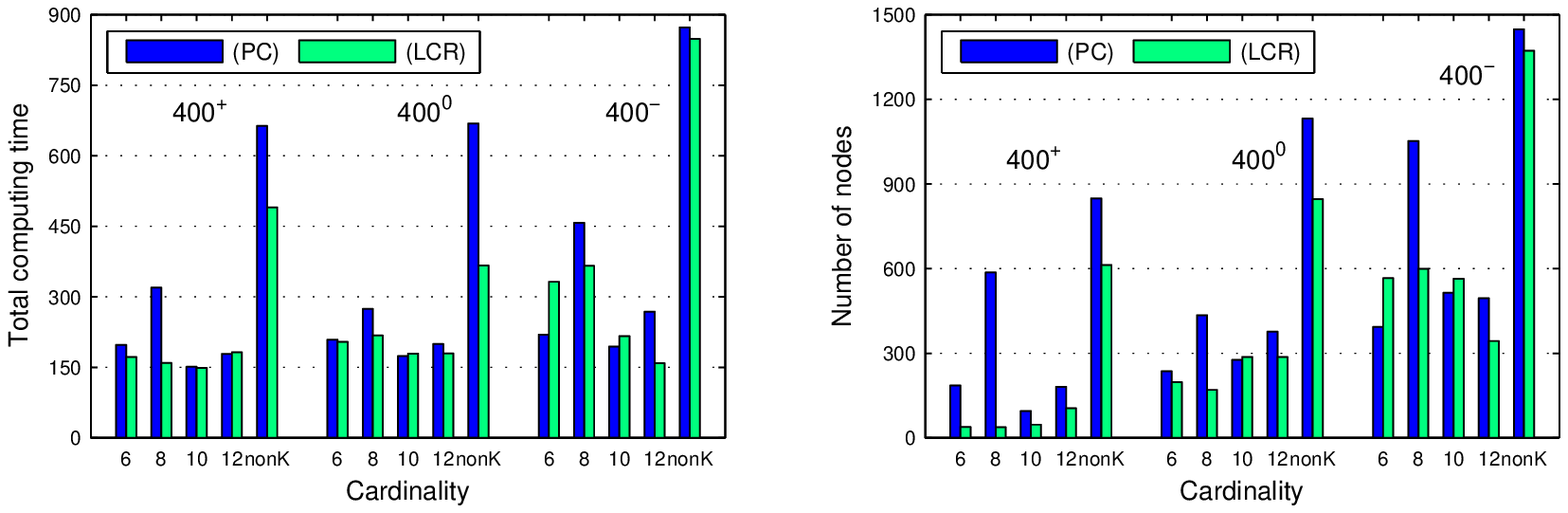}
\caption{\label{mvfigure}Comparison of total computing time and nodes for {\rm (MV)}. }
\end{center}
\end{figure}

Figure \ref{mvfigure} displays the total computing time and nodes of the two
reformulations for (MV).
The total computing time for (PC) is the sum of ``$({\rm time}_l)$'' and the ``time'' for (PC), and
the total computing time for (LCR) is the sum of ``$({\rm time}_l)$'', ``$({\rm time}_s)$''
and the  ``time'' for (LCR).

From Figure \ref{mvfigure}, we can see that, in terms of the total computing time, (LCR) performs better than (PC) for 25 out of the total 45 cases.
If we omit the time for solving the SDP and SOCP and only compare the time of
the MIQP solver \texttt{CPLEX}, (LCR) performs better than (PC) for 41 out of the total  45 cases.
As the perspective cut approach represents the state-of-the-art,
the test result for this (MV) data set confirms that using (LCR) reformulation to solve (P) is also efficient.
We need to emphasize that, we have tried our best in our numerical tests to optimize the implementation of the
perspective cut approach, as the efficiency of the perspective cut approach
depends heavily on its implementation details and also on selected
parameters of the MIQP solvers.

For the number of nodes explored in (LCR) and (PC), they
closely match each other.
We might think that because
the relaxations in (LCR) are in general looser than the ones in (PC) at children nodes,
the number of nodes explored by (LCR) should be larger than (PC)
all the time.
However, this might not always be the case because the branching
schemes, feasible solution heuristics
and the branch-and-bound tree inside \texttt{CPLEX}
could be quite different for (PC) and (LCR)
and the number of nodes explored could demonstrate a more random pattern.

\subsection{Subset selection problem}

In this subsection, we compare (PC) and (LCR) for the subset selection problem (SSP)
introduced in the introduction section.

We use the $40$ instances of the subset selection problem from \cite{zheng2013}
with $n=50,100$ and $K=5,10,15,20$, $5$ instances for each $(n,K)$ pair.
In those instances, we set $m=2n$. The elements of $A$ are generated from the standard
normal distribution $N(0,1)$ and $b=A\beta+\epsilon$ where the elements of $\epsilon$
are generated from the standard normal distribution $N(0,1)$ and the elements of $\beta$
are generated uniformly form $[-1,1]$.
The lower and upper bounds for the solution $x_i,~i=1,\ldots,n$, are set, respectively, at
$-100$ and $100$, which are sufficiently large for those instances.

\begin{table}[t!] \footnotesize \vspace*{0.0in} \centering
\caption{\label{tab2}Numerical results of reformulations for {\rm (SSP)}}
\tabcolsep=8pt
\begin{tabular}{cccccccccccc}
\hline \noalign{\smallskip}
     \multirowcell{3}{$n$}
   & \multirowcell{3}{$K$}
   & \multirowcell{3}{$({\rm time}_l)$}
   & \multirowcell{3}{$({\rm time}_s)$}
   &&\multicolumn{2}{c}{${\rm (PC)}$}
   &&\multicolumn{2}{c}{${\rm (LCR)}$}\\
     \noalign{\smallskip}
     \cline{6-7} \cline{9-10}
     \noalign{\smallskip}
   &&&
   && time & nodes
   && time & nodes \\
\hline \noalign{\smallskip}
   \multirowcell{4}{$50$}
   & 5 & 3.86  & 0.14  &  & 2.90  & 61  &  & 0.35  & 115  \\
   & 10 & 3.58  & 0.21  &  & 9.17  & 205  &  & 1.29  & 445  \\
   & 15 & 3.24  & 0.20  &  & 6.82  & 119  &  & 0.71  & 218  \\
   & 20 & 3.44  & 0.25  &  & 9.20  & 136  &  & 0.77  & 276  \\
\hline \noalign{\smallskip}
   \multirowcell{4}{$100$}
   & 5 & 14.39  & 0.48  &  & 8.69  & 99  &  & 3.35  & 295  \\
   & 10 & 14.90  & 0.37  &  & 34.65  & 334  &  & 14.11  & 993  \\
   & 15 & 14.95  & 0.49  &  & 89.46  & 739  &  & 41.38  & 3767  \\
   & 20 & 13.36  & 0.38  &  & 613.27  & 3737  &  & 140.76  & 11789  \\
\hline \noalign{\smallskip}
\end{tabular}
\end{table}

Table \ref{tab2} summarizes the numerical results for the 40 instances of
(SSP).
Each line reports the average results for $5$ instances for each $(n,K)$ pair.
The notations in the table are the same as those in Table \ref{tab1}.

\begin{figure}[ht!]
\begin{center}
\includegraphics[width=5in]{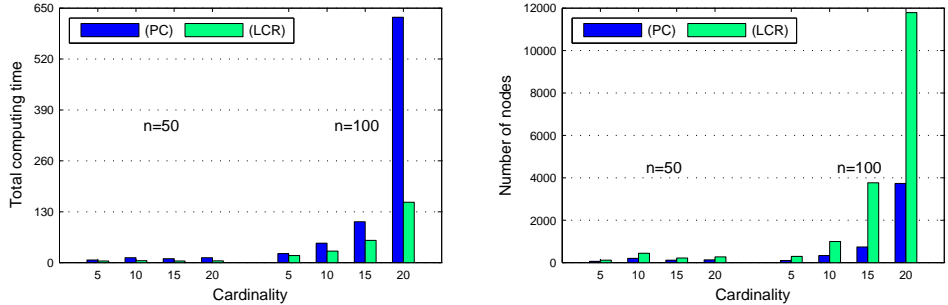}
\caption{\label{sspfigure}Comparison of total computing time and nodes for {\rm (SSP)}.}
\end{center}
\end{figure}

Figure \ref{sspfigure} displays the total computing time and nodes of the two
reformulations for (SSP).
The total computing time for (PC) is the sum of ``$({\rm time}_l)$'' and the ``time'' for (PC),
and the total computing time for (LCR) is the sum of ``$({\rm time}_l)$'', ``$({\rm time}_s)$''
and the  ``time'' for (LCR).

From Figure \ref{sspfigure}, we can see that
although (LCR) explores more nodes than (PC),
the total computing time for (LCR) is
smaller than that of (PC) in all cases and this advantage of (LCR) over (PC) becomes more
apparent as the problem size grows and/or as the cardinality increases.
Here for the (SSP) data set,
trading off tighter children-node bounds with faster  processing time
indeed has a good payoff.
We also remark again that
the efficiency of the perspective cut approach
 depends heavily on its implementation details and also on the selected
parameters of the MIQP solvers.
We, however, believe that our new reformulation
derived from the lift-and-convexification approach provides
a good supplement to the current state-of-the-art approaches.

\section{Combination of lift-and-convexification and QCR}
In this section, we further combine our lift-and-convexification approach and
the QCR approach to derive an even tighter reformulation for (P).

Hammer and Rubin \cite{hammer1970} pioneered the QCR approach in the following binary quadratic programs:
\begin{align}
({\rm BQP})~~
\min
&~~x^TQx+c^Tx\non\\
{\rm s.t.}
&~~Ax=d,\non\\
&~~x_i\in\{0,1\},~i=1,\ldots,n, \non
\end{align}
where $Q$ is indefinite.
In the proposed QCR, they added to the
objective function a term $\sum_iu(x_i^2-x_i)$, where $u$ is a scaler and is chosen to be
the negative value of the smallest eigenvalue of $Q$. Billionnet and Elloumi
\cite{billionnet2007} improved this method by adding
the term $\sum_iu_i(x_i^2-x_i)$ with $u_i$ being the optimal
dual variables of a certain semi-definite program (SDP).
Plateau\cite{plateau2006reformulations} and Billionnet et al.\cite{billionnet2008, billionnet2009qcr}
also utilized the equality $Ax=d$ in QCR and added
the term $\sum_iu_i(x_i^2-x_i) + (Ax-d)^T\dg(w)(Ax-d)$
to the objective, where $u$ and $w$ are chosen to be the dual variables
of an enlarged SDP program.
Ahlat{\c{c}}{\i}o{\u{g}}lu et al. \cite{ahlatcciouglu2012combining} proposed to combine QCR and the convex hull relaxation
to solve problem (BQP). The geometric investigation in Li et al. \cite{li2012} for binary quadratic programs provides some theoretical support for QCR from another angle.
Billionnet et al.\cite{billionnet2012qcrExtended} extended the QCR approach to general mixed-integer programs by
using binary decomposition.

To make our discussion more general, we add equality constraints to (P)
and consider the following variant of (P):
\begin{align*}
({\rm P'})~~~
\min
&~~f(x,y):=x^TQx+c^Tx+h^Ty\non\\
s.t.
&~~Ax+By\le d,\non\\
&~~Ex+Fy=   g,\non\\
&~~a_iy_i\le x_i\le b_i y_i,~y_i\in\{0,1\},~i=1,\ldots,n,
\end{align*}
where $g\in\Re^M$, and $E,F\in\Re^{M\times n}$.

%If we apply QCR directly to the constraints
%$a_iy_i\le x_i\le b_i y_i,~y_i\in\{0,1\},~i=1,\ldots,n,$ in $({\rm P'})$,
%we need to introduce first new variables $r_i\ge 0,t_i\ge0$ to transform these constraints to the following equalities,
%\begin{align*}
%a_iy_i-x_i+r_i=0,\\
%x_i-b_iy_i+t_i=0,\\
%\end{align*}
%We can then add the following term,
%\begin{align}
%\sum_{i=1}^nu_i(y_i^2-y_i)+\sum_{i=1}^nv_i(a_iy_i-x_i+r_i)^2+\sum_{i=1}^nw_i(x_i-b_iy_i+t_i)^2
%\label{qcr_1}
%\end{align}
%to the objective function in $({\rm P'})$.
%However, this application of QCR does not help
%because the optimal parameter $u_i$ would be zero, as negative $u_i$ makes the objective function nonconvex and
%positive $u_i$ would loosen the lower bound due to $y_i^2-y_i\le 0$ for $y_i\in[0,1]$).

QCR would become beneficial when being applied to the equality and inequality constraints in (P') on a top of our lift-and-convexification reformulation.
Let us consider now the following equivalent reformulation of $({\rm P'})$:
\begin{align*}
  ({\rm P}(u,v,w,t))~~
  \min
  &~~f(x,y)+\sum_{i=1}^nq_i(x_i,y_i)+(Ex+Fy-g)^T\dg(w)(Ex+Fy-g)\non\\
  &~~~+(Ax+By+s-d)^T\dg(t)(Ax+By+s-d)\\
  s.t.
  &~~Ax+By+s= d,s\ge 0,\non\\
  &~~Ex+Fy=   g,\non\\
  &~~a_iy_i\le x_i\le b_i y_i,~y_i\in\{0,1\},~i=1,\ldots,n,
\end{align*}
where $q_i(x_i,y_i)$ is defined in (\ref{q_2}).

The best parameter set $(u,v,w,t)$ can be found by solving the following
problem:
\begin{align*}
&({\rm MAX}uvwt)\\
&\max\{v({\rm\ovl{P}}(u,v,w,t))\mid The~objective
~function~of~({\rm P}(u,v,w,t))~is~convex\},
\end{align*}
where $({\rm\ovl{P}}(u,v,w,t))$ is the continuous relaxation of $({\rm P}(u,v,w,t))$.

Using the same technique in the proof for Theorem \ref{thmsdpq},
we can convert the problem $({\rm MAX}uvwt)$ to an SDP problem.
\begin{thm}
\label{thmsdpa}
The problem $({\rm MAX}uvwt)$ is equivalent to the following SDP problem:
\begin{align}
  ({\rm SDP}_a)~~
  \max
  &~~\tau \non \\
  {\rm s.t.}
  &~
  \left(
  \begin{array}{cccc}
  P_{11}   & P_{12}   & P_{13}   & P_{14} \\
  P_{12}^T & P_{22}   & P_{23}   & P_{24} \\
  P_{13}^T & P_{23}^T & P_{33}   & P_{34} \\
  P_{14}^T & P_{24}^T & P_{34}^T & P_{44} \\
  \end{array}
  \right)\succeq 0,\label{dsA31} \\
  &~~(\delta,\mu,\sigma,\lambda,\pi)\in\Re^m_+\times\Re^n_+\times\Re^n_+\times\Re^n_+\times\Re^n_+,\non\\
  &~~(\eta,\zeta,u,v,\tau)\in \Re^m\times \Re^M\times\Re^n\times\Re^n\times \Re, \non
\end{align}
where
\begin{align}
   P_{11}&=Q+E^T\dg(w)E+A^T\dg(t)A \non \\
   P_{12}&=\frac{1}{2}{\rm diag}(u)+E^T\dg(w)F+A^T\dg(t)B\non \\
   P_{13}&=A^T\dg(t)\non \\
   P_{14}&=\frac{1}{2}(c-u+A^T\eta+E^T\zeta-2E^T\dg(w)g+2A^T\dg(t)(-d)-\mu+\sigma)\non \\
   P_{22}&={\rm diag}(v)+F^T\dg(w)F+B^T\dg(t)B\non \\
   P_{23}&=B^T\dg(t)\non \\
   P_{24}&=\frac{1}{2}(h-v+B^T\eta+F^T\zeta-2F^T\dg(w)g+2B^T\dg(t)(-d)\non\\
         &~~~+{\rm diag}(a)\mu-{\rm diag}(b)\sigma-\lambda+\pi)\non \\
   P_{33}&=\dg(t) \non \\
   P_{34}&=\frac{1}{2}(\eta-\delta-2\dg(t)d)\non \\
   P_{44}&=-\eta^T d-e^{T}\pi-\zeta^Tg+g^T\dg(w)g+d^T\dg(t)d \non
\end{align}
\end{thm}
\proof
We first express $({\rm\ovl{P}}(u,v,w,t))$ by its dual form.
Associate the following multipliers to the constraints in $({\rm\ovl{P}}(u,v,w,t))$:
\begin{itemize}
  \item $\eta\in\Re^m$ for $Ax+By+s=d$;
  \item $\delta\in\Re^m_+$ for $s\ge 0$;
  \item $\zeta\in\Re^M$ for $Ex+Fy=g$;
  \item $\mu_i$ and $\sigma_i\in \Re_+$ for
        $a_iy_i\le x_i$ and $x_i\le b_iy_i$, respectively, $i=1,\ldots,n$;
  \item $\lambda_i$ and $\pi_i\in \Re_+$ for $y_i\ge 0$ and $y_i\le 1$, respectively,
        $i=1,\ldots,n$.
\end{itemize}
Let $\mu=(\mu_1,\ldots,\mu_n)^T$, $\sigma=(\sigma_1,\ldots,\sigma_n)^T$,
$\lambda=(\lambda_1,\ldots,\lambda_n)^T$ and $\pi=(\pi_1,\ldots,\pi_n)^T$.
Let $\omega$ denote the vector formed by $\delta,\mu,\sigma,\lambda$ and $\pi$.
The Lagrangian function of $({\rm\ovl{P}}(u,v,w,t))$ is then given by

\begin{align*}
  &L(x,y,\omega,\eta,\zeta)\\
  &= x^TQx+c^Tx+h^Ty+\sum_{i=1}^n(u_ix_iy_i+v_iy_i^2-u_ix_i-v_iy_i)\\
  &~~~+(Ex+Fy-g)^T\dg(w)(Ex+Fy-g)\non\\
  &~~~+(Ax+By+s-d)^T\dg(t)(Ax+By+s-d)\\
  &~~~+\eta^T(Ax+By+s-d)+\delta^T(-s)+\zeta^T(Ex+Fy-g)\\
  &~~~+\sum_{i=1}^{n}\mu_i(a_iy_i-x_i)+\sum_{i=1}^{n}\sigma_i(x_i-b_iy_i)
      +\sum_{i=1}^{n}\lambda_i(-y_i)+\sum_{i=1}^{n}\pi_i(y_i-1)\\
  &=  (x,y,s)^T
      \left( \begin{array}{ccc}
             P_{11} & P_{12} & A^T\dg(t)\\
             P_{12}^T & P_{22}  & B^T\dg(t)\\
             \dg(t)A & \dg(t)B & \dg(t) \\
             \end{array}
      \right)
      \left( \begin{array}{c}
             x \\
             y \\
             s \\
             \end{array}
      \right)\\
  &~~~+(c-u+A^T\eta+E^T\zeta-2E^T\dg(w)g+2A^T\dg(t)(-d)-\mu+\sigma)^{T}x\\
  &~~~+(h-v+B^T\eta+F^T\zeta-2F^T\dg(w)g+2B^T\dg(t)(-d)\non\\
  &~~~+{\rm diag}(a)\mu-{\rm diag}(b)\sigma-\lambda+\pi)^Ty
      +(\eta-\delta-2\dg(t)d)^Ts \\
  &~~~-\eta^T d-e^{T}\pi-\zeta^Tg+g^T\dg(w)g+d^T\dg(t)d.
\end{align*}

Furthermore, the Lagrangian dual problem of $({\rm\ovl{P}}(u,v,w,t))$ can be expressed as
\begin{align}
\max \{\min_{(x,y,s)\in\Re^n\times \Re^n\times \Re^m} ~L(x,y,\omega,\eta,\zeta)\mid \omega\ge 0\}.\label{dsA1}
\end{align}
Introducing an additional variable $\tau$, we can rewrite (\ref{dsA1}) as
\begin{align}
  \max
  &~~\tau\label{dsA21}\\
  {\rm s.t.}
  &~~\min_{(x,y,s)\in\Re^n\times \Re^n\times \Re^m} ~L(x,y,\omega,\eta,\zeta)\ge\tau,\label{dsA22}\\
  &~~\omega\ge 0.\label{dsA23}
\end{align}
We see that the constraint in (\ref{dsA22}) is equivalent to
$L(x,y,\omega,\eta,\zeta)-\tau\ge 0$ for all $x,y$,s, which is further equivalent to
\begin{align}
  L(x/k,y/k,\omega,\eta,\zeta)-\tau\ge 0,
  ~~\forall (x,y,s)\in\Re^n\times \Re^n\times \Re^m,~\forall k\in\Re,~k\not=0.\label{LhomoA1}
\end{align}
Multiplying both sides of (\ref{LhomoA1}) by $k^2$ yields a homogeneous quadratic form
of $(x,y,s,t)$ in the left-hand side of (\ref{LhomoA1}).
Thus, the constraint in (\ref{dsA22}) is equivalent to the semidefinite constraint (\ref{dsA31}).
Consequently, the problem in (\ref{dsA21})-(\ref{dsA23}) can be expressed as
\begin{align}
\max\{\tau\mid {\rm constraint}~(\ref{dsA31}),~\omega\ge 0\}.\label{dsA4}
\end{align}
If the~objective
~function~of~$({\rm P}(u,v,w,t))$~is~convex, by the strong duality of convex quadratic programming
(see, e.g., Proposition 6.5.6 in \cite{bertsekas2003}),
the optimal values of $({\rm\ovl{P}}(u,v,w,t))$ and (\ref{dsA4}) are equal.
Thus,  we have shown that problem $({\rm MAX}uvwt)$ is equivalent to an SDP problem
in the form of $({\rm SDP}_a)$.
\endproof

Although this SDP problem $({\rm SDP}_a)$ is very large and takes time to solve,
it only needs to be solved once to get the new reformulation.
When the original problem $({\rm P'})$ is very difficult to solve,
solving this SDP to get a better reformulation can gain
overall computational advantage.

To test the effectiveness of the new reformulation, we compare
the bounds of $({\rm SDP}_q)$ and $({\rm SDP}_a)$ on
the portfolio selection problem data set introduced in \S 4.
We use the 30 instances with the least diagonal dominance.
For each instance, we impose additional equality constraints by
dividing the stocks into 10 sections and demanding only
one stock to be invested from each section.
Such an additional constraint is a very practical one,  as in real life applications
portfolios are often constructed by choosing investment opportunities from different industries and sections.
We use the following measure for bound improvement,
\[
{\rm impr.} = \frac{v({\rm SDP}_a)~-~v({\rm SDP}_q)}
                   {{\rm optimal\_objective\_value}~-~v({\rm SDP}_q)}.
\]
Table \ref{tab3} shows the bound improvement for the 30 instances, which ranges from
$2.5\%$ to $47.1\%$, resulting an average bound improvement around $20\%$.
This numerical experiment confirms that combining the lift-and-convexification approach and QCR generates a much tighter reformulation on average.

\begin{table}[t!] \footnotesize \vspace*{0.0in} \centering
\caption{\label{tab3}Numerical results of bound improvement from combining lift-and-convexification and QCR}
\tabcolsep=8pt
    \begin{tabular}{|ccc|ccc|ccc|}
    \hline
    n          & inst. & impr.  & n          & inst. & impr.  & n          & inst. & impr.  \\
    \hline
    $200^-$   & 1     & 17.0\% & $300^-$   & 1     & 47.1\% & $400^-$   & 1     & 15.1\% \\
    $200^-$   & 2     & 28.0\% & $300^-$   & 2     & 12.7\% & $400^-$   & 2     & 21.8\% \\
    $200^-$   & 3     & 31.8\% & $300^-$   & 3     & 25.9\% & $400^-$   & 3     & 6.5\% \\
    $200^-$   & 4     & 18.1\% & $300^-$   & 4     & 17.6\% & $400^-$   & 4     & 10.1\% \\
    $200^-$   & 5     & 12.7\% & $300^-$   & 5     & 26.1\% & $400^-$   & 5     & 16.5\% \\
    $200^-$   & 6     & 17.7\% & $300^-$   & 6     & 22.0\% & $400^-$   & 6     & 17.9\% \\
    $200^-$   & 7     & 31.1\% & $300^-$   & 7     & 22.8\% & $400^-$   & 7     & 42.1\% \\
    $200^-$   & 8     & 29.2\% & $300^-$   & 8     & 15.6\% & $400^-$   & 8     & 22.5\% \\
    $200^-$   & 9     & 30.4\% & $300^-$   & 9     & 2.5\%  & $400^-$   & 9     & 11.5\% \\
    $200^-$   & 10    & 33.0\% & $300^-$   & 10    & 30.2\% & $400^-$   & 10    & 20.1\% \\
    \hline
    average &       & 24.9\% &       &       & 22.2\% &       &       & 18.4\% \\
    \hline
    \end{tabular}%
\end{table}

\section{Concluding remarks}
We have developed in this paper the lift-and-\\ convexification approach to construct a
parameterized set of MIQP reformulations for convex quadratic programs with semi-continuous variables.
The primary idea behind this approach is to lift the quadratic term in the objective function from
the $x$-space to the $(x,y)$-space and to convexify the resulting quadratic function of $(x,y)$.
We have proposed an SDP formulation to identify the best MIQP reformulation from among this parameterized
set and have proved that the identified best reformulation  has a continuous relaxation that is as tight as
the continuous relaxation of the well known perspective reformulation.
By revealing the relationship between our new reformulation and the perspective reformulation,
we further reduce the computational effort required to construct our new reformulation and show that
we only need little extra effort to solve an additional SOCP problem when compared to the perspective reformulation. Most importantly, our new reformulation
retains the linearly constrained quadratic programming structure of the original mix-integer problem, which facilitates more effective utilization of commercial mixed integer programming solvers and ensures much faster computational time
at children nodes in the branch-and-bound searching process.
Our preliminary comparison results indicate that the performance
of our new reformulation solved in general MIQP solvers
is, at least, competitive to the state-of-the-art perspective
cut approach in many cases and provides a good
supplement to the state-of-the-art approaches.
We further combine our lift-and-convexification approach
and the quadratic convex reformulation approach in the literature
to obtain an even tighter reformulation.
In a broader picture, the lift-and-convexification approach offers an efficient solution
framework of tight MIQP reformulation which improves the existing literature on the trade-off between the bound quality and computational complexity.

\bibliographystyle{siam}
\bibliography{qp_semi}

\end{document}